\documentclass[12pt]{amsart}
\usepackage{amsmath,amsthm,latexsym,amscd,amsbsy,amssymb,amsfonts,
euscript}

\newcommand{\St}{\mathop{\rm St}\nolimits}

\def\ed{\text{e\,-}\dim}

\newcommand{\cov}{\mathop{\rm cov}\nolimits}
\newcommand{\ord}{\mathop{\rm ord}\nolimits}

\def\wh{\widehat}

\newtheorem{thm}{Theorem}[section]

\newtheorem{lem}[thm]{Lemma}
\newtheorem{pro}[thm]{Proposition}

\theoremstyle{definition}

\theoremstyle{remark}

\chardef\bslash=`\\ 

\makeatletter
\def\verbatim{\interlinepenalty\@M \@verbatim
  \leftskip\@totalleftmargin\advance\leftskip2pc
  \frenchspacing\@vobeyspaces \@xverbatim}
\makeatother
\numberwithin{equation}{section}

\begin{document}


\title[Universal absolute extensors in extension theory]
{Universal absolute extensors in extension theory}

\author{Alex Karasev}
\address{Department of Computer Science and Mathematics,
Nipissing University, 100 College Drive, P.O. Box 5002, North Bay,
ON, P1B 8L7, Canada} \email{alexandk@nipissingu.ca}
\thanks{The authors were partially supported by their NSERC grants.}

\author{Vesko Valov}
\address{Department of Computer Science and Mathematics, Nipissing University,
100 College Drive, P.O. Box 5002, North Bay, ON, P1B 8L7, Canada}
\email{veskov@nipissingu.ca}




\begin{abstract}
Let $L$ be a countable and locally finite CW complex. Suppose that
the class of all metrizable compacta of extension dimension $\le
[L]$ contains a universal element which is an absolute extensor in
dimension $[L]$. Our main result shows that $L$ is quasi-finite.
\end{abstract}

\maketitle

\markboth{A.~Karasev and V.~Valov}{Universal absolute extensors in
extension theory}


\section{Introduction}

In this note we deal with one of the central problems in
extension theory that can be described as follows. Consider a CW
complex $L$. Suppose that the class of all metrizable compacta of
extension dimension $\le [L]$ has a universal element which is an
absolute extensor in dimension $[L]$. What can be said about the
properties of the complex $L$? It is known \cite[Theorem
2.5]{ch01} that the situation described above occurs when $L$ is a
finite complex. The main purpose of this note is to show that such
a complex $L$ must be necessarily quasi-finite. We do not know
whether this condition is also sufficient. Note that quasi-finite
CW complexes were introduced in \cite{kar} as complexes which
provide the solution of the following problem: characterize all
complexes $P$ such that there exists a $P$-invertible mapping of
metrizable compactum of extension dimension $\le [P]$ onto the
Hilbert cube. Note also that existence of such a mapping for a
complex $P$ implies the existence of a universal metrizable
compactum of extension dimension $[P]$. Consequently, if $L$ is
quasi-finite, it guarantees the existence of universal compactum
of extension dimension $[L]$.

As an application of our result we show that there is no universal
compactum of a given cohomological dimension which is an absolute
extensor with respect to spaces of given cohomological dimension.
Similar result \cite{zar} was known only for the case of integral
cohomological dimension.

\section{Preliminaries}

In this and all subsequent sections ``complex" means ``a countable
and locally finite CW complex". All spaces under consideration are
assumed to be Tychonov and all maps are continuous. The letter
``$L$" will be reserved to denote a complex. In this note, for
spaces $X$ and $Y$, the notation $Y \in \rm AE (X)$ will always
mean that every map $f\colon A \to Y$, defined on a closed
subspace $A$ of $X$, admits an extension over $X$. By $[L]$ we
denote the extension type of a complex $L$ and by $\ed X$ we
denote the extension dimension of space $X$ \cite{D,DD}. For a
normal space $X$, inequality $\ed X\le [L]$ means that $L\in
AE(X)$. More information about extension dimension and extension
types can be found in \cite{ch97,ch01}.

We say that a map $f\colon X\to Y$ is {\it $[L]$-soft} \cite{ch01}
if for each Polish space $Z$ with $\ed Z\le [L]$, for each closed
subspace $A$ of  $Z$, and for any two maps $g\colon Z\to Y$ and
$h\colon A\to X$ such that $f\circ h =g|_A$ there exists a map
$\overline{h}\colon Z\to X$ extending $h$ and satisfying the
conditions $f\circ\overline{h} = g$.

Let $\mathcal B$ be a certain class of spaces. We shall say that a
space $X$ is an absolute extensor in dimension $[L]$ for the class
$\mathcal B$ (notation $X\in AE ([L],{\mathcal B})$) if $X\in AE
(Y)$ for every $Y$ from $\mathcal B$ such that $\ed Y \le [L]$. We
shall denote the class of all metrizable compacta by $\mathcal C$
and the class of all Polish spaces by
$\mathcal P$. 
The following remark is trivial.

\begin{pro}\label{soft}
Let $f\colon X\to Y$ be an $[L]$-soft mapping. Then $X\in AE([L],
\mathcal P)$ iff $Y\in AE([L],\mathcal P)$.
\end{pro}

\noindent Let $X$ be a normal space. A pair of spaces $V\subset U$
is called {\it $X$-connected} if for every closed subspace
$A\subset X$ any mapping of $A$ to $V$ can be extended to a
mapping of $X$ into $U$. Suppose that $\mathcal B$ is a certain
subclass of the class of normal spaces. A pair of spaces $V\subset
U$ is called {\it $[L]$-connected with respect to $\mathcal B$} if
for every space $X\in\mathcal B$ with $\ed X\le [L]$ the pair
$V\subset U$ is $X$-connected. In what follows we will need the
following observation from \cite[Proposition A.1]{brchka} (recall
that in this note we consider only countable complexes).

\begin{pro}\label{lconnected}
Let $L$ be a complex and $V\subset U$ be a pair of Polish spaces.
If this pair is $[L]$-connected with respect to Polish spaces then
it is $[L]$-connected with respect to all normal spaces.
\end{pro}

We say \cite{kar} that a complex $L$ is {\it quasi-finite} if for
every finite subcomplex $P$ of $L$ there exists a finite
subcomplex $P'$ of $L$ containing $P$ such that the pair $P\subset
P'$ is $[L]$-connected with respect to Polish spaces.

The following theorem provides a characterization of quasi-finite
complexes.  Note that equivalences from (a) through (e) were
obtained by Chigogidze in \cite[Theorem~2.1]{ch99} and the
equivalence of these properties to (f) follows from
\cite[Theorem~3.1]{kar}.

\begin{thm}\label{equivalence} Let $L$ be a complex. Then
the following statements are equivalent:

\begin{itemize}
    \item[(a)] $\ed \beta X \le [L]$ whenever $X$ is a space with $\ed X \le [L]$.
    \item[(b)] $\ed \beta X\le [L]$ whenever $X$ is a normal space with $\ed X\le [L]$.
    \item[(c)] $\ed \beta (\oplus\{ X_t \, |\, t\in T\} )\le [L]$ whenever $T$ is an arbitrary indexing set
    and $X_t$, $t\in T$, is a separable metrizable space with $\ed X_t\le [L]$.
    \item[(d)] $\ed \beta (\oplus\{ X_t \, |\, t\in T\} )\le [L]$ whenever $T$ is an arbitrary indexing set
    and $X_t$, $t\in T$, is a Polish space with $\ed X_t\le [L]$.
    \item[(e)] There exists a $[L]$-invertible map $f\colon X\to
I^{\omega}$ where $X$ is a metrizable compactum with $\ed X \le
[L]$.
    \item[(f)] $L$ is quasi-finite.

\end{itemize}

\end{thm}

\section{Results}

Let $\mathcal B$ be a subclass of the class of normal spaces. We
say that a complex $L$ possesses {\it connected pairs property
with respect to $\mathcal B$} if for any metrizable compactum $K$
with $\ed K\le [L]$ there exists a metrizable compactum $C$
containing $K$ such that $\ed C\le [L]$ and the pair $K\subset C$
is $[L]$-connected with respect to $\mathcal B$.

\begin{lem}\label{absextforbeta} Let $T$ be an arbitrary indexing set and
$\{ X_t \, |\, t\in T\}$ be a collection of Polish spaces such
that $\ed X_t\le [L]$ for each $t\in T$. Let $X=\oplus\{ X_t\, |\,
t\in T\} $. Suppose that $K\subset C$ is a pair of metrizable
compacta such that $\ed K\le [L]$. If the pair $K\subset C$ is
$[L]$-connected with respect to Polish spaces then this pair is
$\beta X$-connected.
\end{lem}

\begin{proof}
Let $A$ be a closed subset of $\beta X$ and $f\colon A\to K$ be a
map. Consider the adjunction space $Y = X\cup_f f(A)$. Note that
$Y$ can be viewed as the disjoint union of two subspaces,
homeomorphic to $(X-A)$ and $f(A)$, respectively. We claim that
$\ed Y\le [L]$. Indeed, $f(A)$ is a closed subspace of $K$ and
therefore $\ed f(A)\le [L]$. Further, $X - A$ is an open subset of
$X$. Note that $X$ is metrizable and therefore perfectly normal.
Therefore the claim follows from the countable sum theorem.
Observe also that Proposition \ref{lconnected} allows us to assume
that the pair $K\subset C$ is $L$-connected with respect to normal
spaces. Hence the identity mapping $i$ of a copy of $f(A)$ in $Y$
to a copy of $f(A)$ in $K$ can be extended to a mapping $j\colon
Y\to C$. Let $p\colon X\cup A\to Y= X\cup_f f(A)$ be the natural
projection and let $\overline{f} = j\circ p$. Then
$\overline{f}\colon X\cup A\to C$ extends $f$ to $X\cup A$. Now
the unique extension $\beta \overline{f}$ of $\overline{f}$ over
$\beta X$ yields the required extension of $f$.
\end{proof}

\begin{lem}\label{extensors} Let $L$ be a complex possessing the connected
pairs property with respect to Polish spaces and $X$ be a
compactum. Suppose that each pair $K\subset C$ of metrizable compacta with $\ed C\le [L]$ is
$X$-connected provided $K\subset C$ is $[L]$-connected with respect to Polish spaces. Then for every metrizable space $Y\in AE([L],\mathcal C)$ with $\ed Y\le[L]$  we have  $Y\in
AE(X)$
\end{lem}

\begin{proof}
Let $A$ be a closed subset of $X$ and $f\colon A\to Y$ be a
mapping. Note that $f(A)$ is a metrizable compactum and $\ed
f(A)\le [L]$. Therefore there exists a metrizable compactum $B$
with $\ed B\le [L]$ such that the pair $f(A)\subset B$ is
$[L]$-connected with respect to Polish spaces. Hence, by our hypotheses,
the pair $f(A)\subset B$ is $X$-connected. Because $f$ can
be viewed as a map sending $A$ to a copy of $f(A)$ inside $B$,
this map can be extended to a map
$f'\colon X\to B$. Since $Y\in AE([L],\mathcal C)$ the
homeomorphism identifying a copy of $f(A)$ in $B$ with a copy of
$f(A)$ in $Y$ can be extended to a mapping $h\colon B\to Y$.
Clearly the map $\overline{f} = h\circ f'\colon X\to Y$ is an
extension of $f$.
\end{proof}

Everywhere below by $\cov (X)$ we denote the set of all open
covers of a space $X$.
If $A$ is a subset of $X$ and $\omega\in\cov (X)$ we denote the
star of $A$ with respect
to $\omega$ by $\St (A,\omega )$. We say that $\nu\in\cov
(X)$ is a strong star-refinement of  if for
each $V\in\nu$ there exists $W\in\omega$ such that $\St
(V,\nu)\subset W$. The following set of notations is borrowed from
\cite{brchka}. For a cover $\Sigma\in\cov (X)$ we denote by
$\Sigma^{(k)}$ its ``$k$-dimensional skeleton", i.e. the set of
all points in $X$ at which order of $\Sigma$ is at most $k+1$.
Thus we let $\Sigma^{(k)} =\{x\in X\mid\ord_\Sigma x\le k+1\}$.
For elements $s_0,s_1,\dots,s_n\in\Sigma$ with non-empty
intersection $\cap_{i=0}^n s_i$ we define a ``closed
$n$-dimensional simplex"
\[ [s_0,s_1,\dots,s_n]=\bigcup_{i=0}^n s_i\setminus
 \bigcup \{ s\in\Sigma \, |\, s \ne s_i , i=0,1,\dots , n\}\]
and its ``interior" $\langle s_0,s_1,\dots,s_n\rangle=\cap_{i=0}^n
s_i\cap\Sigma^{(n)}$. It is easy to check that the $n$-skeleton
consists of $n$-simplices
\[\Sigma^{(n)}=\bigcup\{[s_{i_0},s_{i_1},\dots,s_{i_n}]\mid
 \cap_{k=0}^n s_{i_k}\ne\varnothing\}\]
and that any ``simplex" consists of its ``boundary" and its
``interior"
\[ [s_0,s_1,\dots,s_n]=\bigcup_{m=0}^n [s_0,\dots,\wh s_m,\dots, s_n]
 \cup\langle s_0,s_1,\dots,s_n\rangle.\]
Clearly, $\Sigma^{(k)}$ is closed in $X$ and $\Sigma^{(n)}=X$ if
the cover $\Sigma$ has order $n+1$. Note also that the
``interiors" of distinct $k$-dimensional ``simplices" are mutually
disjoint and
$$\Sigma^{(k)}=\bigcup\{\langle s_{i_0},s_{i_1},\dots,s_{i_n}\rangle\mid
 \cap_{k=0}^n s_{i_k}\ne\varnothing\}\cup\Sigma^{(k-1)}$$ 

The following lemma can be interpreted as a ``weak" version of
Lemma 3.10 from \cite{karval}.

\begin{lem}\label{lifting} Let $X$ be a compactum and $Z$ be a paracompact space such
that any compact subspace of $Z$ is finitely-dimensional in the
sense of usual Lebesgue dimension. Let $g\colon Y\to Z$ be a
surjection with the following property: for every $z\in Z$ and its
neighborhood $U(z)$ in $Z$ there exists a smaller neighborhood
$V(z)$ of $z$ such that $g^{-1}(V(z))\in AE(X)$. Then for any
$\omega\in\cov (Z)$ and for any mapping $f\colon X\to Z$ there
exists a map $\tilde{f}\colon X\to Y$ such that the maps $f$ and
$g\circ\tilde{f}$ are $\omega$-close.
\end{lem}

\begin{proof}
Note that $f(X)\subset Z$ is compact and therefore $\dim f(X) = n
<\infty$ for some $n$. We let $\omega _0 =\omega$ and inductively
construct a sequence of covers $\omega _1, \omega_2, \dots ,\omega
_n$ as follows. Suppose $\omega _i\in\cov (Z)$ is already
constructed and let $\nu$ be a strong star-refinement of $\omega
_i$.  For each $z\in Z$ we choose $U(z)\in\nu$ containing $z$ and
find a smaller neighborhood $V(z)\subset U(z)$ of $z$ having the
property

$$g^{-1}(V(z))\in AE(X)\eqno{(\dag )}$$

\noindent We let $\omega _{i+1} = \{V(z)\, |\, z\in Z\}$.
Obviously, $\omega _{i+1}$ is a strong star-refinement of $\omega
_{i}$.

Let $\Sigma\in\cov (f(X))$ be a finite strong star-refinement of
$\omega _{n}$ restricted on $f(X)$ such that $\ord\Sigma\le n+1$.
We put $\widehat{\Sigma} = \{f^{-1} (U)\, | \, U\in\Sigma\}$.
Clearly $\widehat{\Sigma}$ is a finite open cover of $X$ of order $\le
n+1$. By induction we construct a sequence of maps $f_0, f_1,
\dots ,f_n$ such that $f_k\colon\widehat{\Sigma} ^{(k)}\to Y$ with
the property

$$g(f_k (x))\in \St (f(x), \omega _{n-k})\eqno{(*)}$$

\noindent for all $k$. In order to construct $f_0$ for each
element $s\in\widehat{\Sigma}$ we choose a point $P_s\in g^{-1}
(f(s))$ and then for every ``closed one-dimensional simplex" $[s]$
we let $f_0 | _{[s]} = P_s$. Suppose that $f_k$ has already  been
constructed. It suffices to define $f_{k+1}$ on the ``interior"
$\langle\sigma\rangle$ of each ``simplex" $[\sigma ]=[
s_0,s_1,\dots, s_{k+1} ]$. Since $\widehat{\Sigma}$ is finite and
the ``interiors" of ``closed $k$-dimensional simplices" are
mutually disjoint we can consider each simplex independently. Let
$[\sigma]' = [\sigma]\cap\widehat{\Sigma} ^{(k)}$. Since $\Sigma$
is a strong star-refinement of $\omega _{n}$ (and consequently of
$\omega _{n-k}$) and because of property $(*)$, we can find
$W_{\sigma}\in\omega _{n-k}$ such that $g(f_k
([\sigma]'))\subset\St (W_{\sigma}, \omega _{n-k})$. Since $\omega
_{n-k}$ is a strong star-refinement  of $\omega _{n-k-1}$ and by
the construction of $\omega _{n-k-1}$ there exists an element
$V_{\sigma}\in\omega _{n-k-1}$ possessing property $(\dag )$ and
such that $\St (W_{\sigma}, \omega _{n-k})\subset V_{\sigma}$.
Since $g^{-1} ( V_{\sigma} )\in AE (X)$ we can extend the mapping
$f_k| _{[\sigma ]'}$ to a mapping $f_{k+1} \colon [\sigma] \to
g^{-1} ( V_{\sigma} ) \subset Y$. It is easy to check that
property $(*)$ is satisfied for $f_{k+1}$.

Finally, we let $\widetilde{f} = f_n$.

\end{proof}

\noindent The following theorem provides a characterization of
quasi-finite complexes in terms of connected pairs property.

\begin{thm}\label{cpproperty} A complex $L$ possesses the connected
pairs property with respect to Polish spaces iff $L$ is
quasi-finite.
\end{thm}
\begin{proof} The ``if" part follows from \cite[Proposition 2.4]{karval}.
In order to establish the ``only if" part we shall show that $L$
satisfies property (d) from Theorem \ref{equivalence}. Let $\{ X_t
\, |\, t\in T\}$ be a collection of Polish spaces where $T$ is an
arbitrary indexing set and assume that $\ed X_t\le [L]$ for each
$t\in T$. Let $X=\oplus\{ X_t\, |\, t\in T\} $. We need to show
that $\ed \beta X\le [L]$. Let $A$ be a closed subset of $\beta X$
and $f\colon A\to L$ be a map. Consider an $[L]$-soft mapping
$g\colon Y\to L$ where $Y$ is a Polish space with $\ed Y\le [L]$.
Existence of such mapping follows from \cite[Proposition
5.9]{ch97}. We claim that the mapping $g$ satisfies conditions of
Lemma \ref{lifting}.
Indeed, consider $z\in L$ and its open neighborhood $U(z)$. Then
$U(z)$ contains a neighborhood $V(z)$ of $z$ in $L$ which is an
absolute extensor. Propositions \ref{soft} implies that
$g^{-1}(V(z))\in AE([L],\mathcal P)$. Subsequently applying Lemma
\ref{absextforbeta} and Lemma \ref{extensors} (for the pair
$g^{-1}(V(z))\subset g^{-1}(V(z))$), we conclude that
$g^{-1}(V(z))\in AE(\beta X)$. This proves the claim. Note also
that the same arguments show that $Y\in AE(\beta X)$.

Since $L$ is $ANR$-space there exists an open cover $\omega\in\cov
(L)$ such that any two $\omega$-close maps to $L$ are homotopic.
Applying Lemma \ref{lifting} to mappings $g\colon Y\to L$,
$f\colon A\to L$, and to the cover $\omega$ we obtain a map
$\tilde{f}\colon A\to Y$ such that $f$ and $g\circ\tilde{f}$ are
$\omega$-close. Since $Y\in AE(\beta X)$ we can extend $\tilde f$
to a map $\overline{f}\colon \beta X\to Y$. Let $f' = g\circ
\overline{f}\colon\beta X\to L$. Note that $f'|_A$ is
$\omega$-close to $f$ and therefore $f$ admits an extension over
$\beta X$, as required.
\end{proof}

For a given complex $L$ we let $\mathcal C _{L}$ to be the class
of all metrizable compacta of extension dimension $\le [L]$. We
say that $X_{L}$ is a universal element for $\mathcal C_{L}$ if
$X_L$ is a metrizable compactum with $\ed X_L\le [L]$ which
contains a topological copy of any metrizable compactum of
extension dimension $\le [L]$. The theorem below contains the main
result of this note and follows directly from Theorem
\ref{cpproperty}.

\begin{thm}\label{main} Let $L$ be a complex and
$\mathcal C_{L}$ be the class of all metrizable compacta of
extension dimension $\le [L]$. If $\mathcal C_{L}$ contains a
universal element $X_{L}$ with the property $X_{L}\in
AE([L],\mathcal P)$ then $L$ is quasi-finite.
\end{thm}

It follows from \cite[Corollary~2.2]{karval} that none of the
Eilenberg-MacLane complexes $K(G,n)$, $n\geq 2$ and $G$ an Abelian
group, is quasi-finite. Therefore Theorem \ref{main} implies the
following result.

\begin{thm} Let $G$ be a countably generated abelian group and $n$ be an integer,
$n\ge 2$. There is no universal compactum of given cohomological
dimension $n$ with respect to the coefficient group $G$, which is
an absolute extensor with respect to Polish spaces of
cohomological dimension $\le n$.
\end{thm}

\noindent In the case of integral cohomological dimension this
theorem is  similar to the observation, made by Zarichnyi in
\cite{zar}.




\end{document}